\documentclass[oneside,letterpaper,spanish,activeacute,english,links]{amca}
\usepackage{graphicx}
\usepackage[latin1]{inputenc}
\usepackage{amsmath,amstext,amsfonts,amssymb,amsbsy,amscd,mathrsfs,mathrsfs,array,multirow} 
\usepackage{pstricks}

\newcommand{\DD}{\mathbf{D}}
\newcommand{\GG}{\mathbf{G}}
\newcommand{\BB}{\mathbf{B}}
\newcommand{\uu}[1][H]{\ensuremath{u_{_\mathrm{#1}}}}
\newcommand{\VV}{\mathscr{V}}
\newcommand{\VVH}[1][H]{\mathscr{V}_{\mathrm{#1}}}
\newcommand{\dOmega}{{\rm{d{\scriptstyle\Omega}}}}

\decimalpoint

\title{UN ANÁLISIS COMPARATIVO DE LOS MÉTODOS MIMÉTICOS, DIFERENCIAS FINITAS Y ELEMENTOS FINITOS PARA PROBLEMAS ESTACIONARIOS}
\author[a,b,1]{Abdul Lugo}
\author[b]{Giovanni Calderón}
\affil[a]{Universidad Polit\'ecnica Territorial del Oeste de Sucre ``Clodosbaldo Russian'', Cuman\'a 6101, Venezuela}
 \affil[b]{Grupo Ciencias de la Computación, Departamento de Matemáticas, Facultad de Ciencias, Universidad de Los Andes, La Hechicera, Mérida 5101, Venezuela}

\begin{document}

\maketitle
\selectlanguage{spanish}

\noindent\hrulefill

\begin{abstract}
Los métodos numéricos: miméticos, diferencias finitas y elementos finitos, son analizados desde un punto de vista numérico. Se busca concluir sobre la eficiencia, orden de convergencia y costo computacional de estos métodos.  El análisis es hecho en  problemas de valor de frontera unidimensionales (ecuación convección-difusión en régimen estacionario) con variaciones distintas en el gradiente, coeficiente difusivo y velocidad convectiva.

\noindent{\bf Palabras clave:}
Métodos miméticos, Método de los elementos finitos, Método de diferencias finitas, Métodos conservativos, Convergencia.

{\footnotetext[1]{Autor para Correspondencia: \url{abdull@ula.ve}}}
\end{abstract}

\vspace{.6cm}

\selectlanguage{english}
\begin{abstract}
Numerical methods: mimetic finite differences and finite elements, are analyzed from a numerical point of view. It seeks to conclude on the efficiency, order of convergence and computational cost of these methods. The analysis is done in boundary value problems one-dimensional (convection-diffusion equation at steady) with different variations in the gradient, diffusion coefficient and convective velocity.


\end{abstract}

\selectlanguage{spanish}

\section{Introducción}
Los métodos numéricos comúnmente usados para definir aproximaciones numéricas a la solución de los diversos problemas de la ingeniería y las ciencias, que se modelan a partir de ecuaciones diferenciales parciales son: el método de diferencias finitas (DF) y el método de los elementos finitos (MEF). No obstante, en las últimas dos décadas, un nuevo tipo de esquemas conservativos en diferencias finitas, conocido originalmente como operadores de soporte y posteriormente como métodos miméticos (MIM) \citep{Shashkov1994,Hyman97,Hyman2002,Castillo2003}, a mostrado su superioridad ante los esquemas clásicos de diferencias finitas \citep{Freites2004,Guevara2005}.

Cada uno de estos métodos presenta ventajas para algún tipo de problema y desventajas en otros. A pesar de esto, y sin introducir ningún tipo de sesgo, el MEF ha predominado en casi todas las aplicaciones que resultan de interés científico hoy en día. Y, en aquellas donde el método DF mantuvo su preponderancia, el MEF estándar ha sufrido modificaciones para optimizar su aplicabilidad en estos problemas. Por ejemplo, se puede citar a \citep{Li2006,Beatrice2008} para problemas de dinámica de fluidos y \citep{Cordero2010} en el caso de problemas con conductividad discontinua.   Por otro lado, mucho se ha dicho sobre la eficiencia de los métodos miméticos en comparación con el método DF \citep{Guevara2008,Guevara2005b}, principalmente en problemas difusivos estacionarios. Sin embargo, se ha dejado de lado su comparación con el MEF, justificando el hecho a una alta complejidad de la implementación computacional del mismo.

En este trabajo, buscando una mejor interpretación de las diferencias que existen entre los tres métodos, se presenta un análisis numérico de las propiedades que posee cada uno de estos métodos: convergencia, precisión en la frontera del dominio, eficiencia en diversos problemas, flexibilidad al dominio físico, implementación y costo computacional. Dicho análisis se enmarca dentro de la ecuación de convección-difusión estacionaria para problemas unidimensionales; se busca facilitar el análisis de los métodos, dejando de lado, momentáneamente, las dificultades analíticas y geométricas que implica un dominio físico 2D.

El resto del artículo se estructura de la siguiente manera. En el siguiente apartado se introduce el problema de valor de frontera usado junto a sus condiciones de frontera (Robin). En el tercer apartado se describen brevemente los métodos numéricos; posteriormente, en la cuarta sección, se realiaza la experimentación numérica y la discusión de los resultados. Finalmente, se presentan las conclusiones obtenidas y las referencias del trabajo.

\section{Problema Modelo}
La ecuación que modela la transferencia de calor por convección-difusión en régimen estacionario en su forma más simple (problema unidimensional
) viene dada por
\begin{equation}\label{ec:1}
    k(x)\frac{d^2}{dx^2}u+\nu(x)\frac{d}{dx}u=f,\qquad\mbox{en }\quad \Omega=(a,b),
\end{equation}
donde $u(x)$ representa la temperatura  (variable del problema) en un punto $x$ del dominio $\Omega=(a,b)$; $k>0$, el coeficiente de difusión térmica;  $\nu$, la velocidad  convectiva (o advectiva) y $f$, una función escalar que describe la existencia de una fuente o sumidero en el problema. La ecuación se completa al definir las condiciones de frontera:
\begin{equation}\label{ec:2}
    \alpha_au(a)+\beta_a\dfrac{d}{dx}u(a)=\gamma_a\qquad\qquad
    \alpha_bu(b)+\beta_b\dfrac{d}{dx}u(b)=\gamma_b
\end{equation}
con $\{\alpha_i,\beta_i,\gamma_i\}$, $i=a,b$, parámetros reales conocidos, y dependiendo de su valor se tendrán condiciones de contorno Dirichlet, Neumann o Robin. Aunque la solución general del problema de valor de frontera \eqref{ec:1}-\eqref{ec:2} se puede obtener en gran parte de las configuraciones del problema mediante métodos analíticos, el mismo resultará apropiado para dejar ver las propiedades de los métodos numéricos en estudio, y simplificar la visualización y análisis de los resultados.

\section{Métodos Numéricos}
En este apartado se describen brevemente los tres métodos utilizados. Una descripción detallada del esquema mimético
a ser utilizado en este artículo puede ser vista, entre otras, en \citep{Guevara2005b}. Para el caso del MEF existe un sin número de referencias que introducen el método desde distintas vertientes; aquí, nos inclinamos por \citep{BeckerCareyOden1981,Solin2006,Calderon2011}. Por último, en el caso de DF se remite al lector a \citep{Strikwerda2004}.

\subsection{Método Mimético}
Los métodos miméticos se basan en la discretización de los operadores clásicos de las EDP (divergencia, gradiente y rotacional) de tal forma que ellos satisfagan una versión discreta del Teorema de Stokes o identidad de Green:
\begin{equation}\label{ec:3}
    \langle \mathbf{D} v,f \rangle_{Q}+\langle v,\mathbf{G} f\rangle_{P}=\langle \mathbf{B }v,f\rangle_{I}.
\end{equation}
Aquí, $\DD$, $\GG$ y $\BB$ son las versiones discretas de sus continuos correspondientes:
gradiente $(\nabla)$, divergencia $(\nabla\cdot)$ y operador de frontera ${\partial}/{\partial n}$. Los $\langle~\rangle$ representan un producto interior generalizado con pesos $Q$, $P$ y $I$.
Usando la identidad \eqref{ec:3} se obtiene una relación para el operador de frontera
\begin{equation}\label{ec:4}
    \BB=Q\DD+\GG^{\mathrm{t}}P.
\end{equation}

\begin{figure}[!h]
 \centering
 \includegraphics[clip,angle=0,width=.7\hsize]{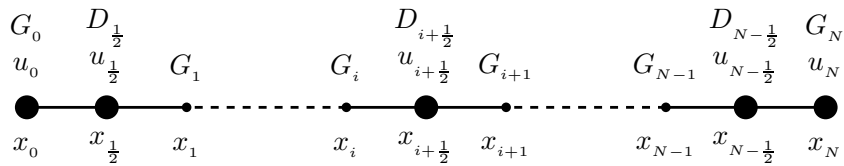}
\caption{Malla unidimensional mimética.}\label{figura:mallaMimetica}
\end{figure}

Para la discretización, se define una malla cuya geometría está dada por los nodos $x_i$, con $i=0,1,\ldots,N$,
y las celdas de la misma son los intervalos $[x_{i-1},x_i]$, para $i=1,\ldots,N$. El tamaño de la celda, $h$, viene dado por $h=1/N$, suponiendo que la malla está distribuida uniformemente. Los nodos intermedios de las celdas queda dado por $x_{i+{1}/{2}}=(x_{i}+x_{i+1})/{2}$,  ver Figura \ref{figura:mallaMimetica}. La solución y el operador divergencia se definen en el centro de las celdas, mientras el operador gradiente en los nodos $x_i$ que definen las celdas (ver Figura \ref{figura:mallaMimetica}).

En este trabajo, se analiza el método mimético que surge de los operadores discretos de segundo orden (tanto en frontera como en los nodos internos) introducidos por \citep{Castillo2003} y dados en \eqref{ec:5} y \eqref{ec:6}.
\begin{equation}\label{ec:5}
\overrightarrow{G}\vec{u}\equiv\left[
                                 \begin{array}{c}
                                 (Gu)_{_{0}} \\
                                 (Gu)_{_{1}} \\
                                 \vdots \\
                                 (Gu)_{_{N-1}} \\
                                 (Gu)_{_{N}} \\
                                 \end{array}
                               \right]=\frac{1}{h}\left[
                                                    \begin{array}{cccccc}
                                                      -\frac{8}{3} & 3 & -\frac{1}{3} & 0 & \cdots & 0 \\
                                                      0 & -1 & 1 & 0 & \cdots & 0 \\
                                                      \vdots & \ddots & \ddots & \ddots & \ddots & \vdots \\
                                                      0 & \cdots & 0 & -1 & 1 & 0 \\
                                                      0 & \cdots & 0 & \frac{1}{3} & -3 & \frac{8}{3} \\
                                                    \end{array}
                                                  \right]_{_{N+1\times N+2}}\left[
                                                                              \begin{array}{c}
                                                                                u_{_{0}} \\
                                                                                u_{_{\frac{1}{2}}} \\
                                                                                \vdots \\
                                                                                u_{_{N-\frac{1}{2}}} \\
                                                                                u_{_{N}} \\
                                                                              \end{array}
                                                                            \right]
\end{equation}

\begin{equation}\label{ec:6}
\overrightarrow{D}\vec{v}\equiv\left[
                                 \begin{array}{c}
                                 0 \\
                                 D_{_{\frac{1}{2}}}v \\
                                 D_{_{\frac{3}{2}}}v \\
                                 \vdots \\
                                 D_{_{N-\frac{3}{2}}}v \\
                                 D_{_{N-\frac{1}{2}}}v \\
                                 0 \\
                                 \end{array}
                               \right]=\frac{1}{h}\left[
                                                    \begin{array}{cccccc}
                                                      0 & 0 & 0 & 0 & \cdots & 0 \\
                                                     -1 & 1 & 0 & 0 & \cdots & 0 \\
                                                      0 &-1 & 1 & 0 & \cdots & 0 \\
                                                      \vdots & \ddots & \ddots & \ddots & \ddots & \vdots \\
                                                      0 & \cdots & 0 & -1& 1 & 0 \\
                                                      0 & \cdots & 0 & 0 & -1& 1 \\
                                                      0 & \cdots & 0 & 0 & 0 & 0 \\
                                                    \end{array}
                                                  \right]_{_{N+2\times N+1}}\left[
                                                                              \begin{array}{c}
                                                                              v_{_{0}} \\
                                                                              v_{_{1}} \\
                                                                              v_{_{2}} \\
                                                                              \vdots \\
                                                                              v_{_{N-2}} \\
                                                                              v_{_{N-1}} \\
                                                                              v_{_{N}} \\
                                                                              \end{array}
                                                                            \right]
\end{equation}

En \eqref{ec:6} aparecen dos filas de ceros (una en la parte superior y otra en la parte inferior), con esto se busca obtener
una matriz cuadrada al componer operadores. 
El operador de frontera $B$ queda dado explícitamente por
\[
 B=\left[\begin{array}{ccccccc}
      -1 & 0    & 0 & \cdots & 0  & 0 & 0 \\
      1/8& -1/8 & 0 & \cdots & 0  & 0 & 0 \\
     -1/8& 1/8  & 0 & \cdots & 0  & 0 & 0 \\
      0  & 0    & 0 & \cdots & 0  & 0 & 0 \\
 \vdots  &\vdots        & \vdots       &\ddots        & \vdots      & \vdots &\vdots  \\
      0  & 0    & 0 & \cdots & 0  & 0    & 0 \\
      0  & 0    & 0 & \cdots & 0  & -1/8 &1/8 \\
      0  & 0    & 0 & \cdots & 0  & 1/8  &-1/8 \\
      0  & 0    & 0 & \cdots & 0  & 0    &1
    \end{array}\right]_{N+2\times N+2}\]

A partir de las discretizaciones de los operadores, la aproximación mimética para la ecuación de
convección-difusión \eqref{ec:1} queda dada por
\begin{equation}\label{ec:7}
    k\DD\GG U+\nu\GG U = (k\DD\GG+\nu\GG)U=F,
\end{equation}
donde $U$  representa la solución aproximada mimética  a la solución exacta, $u$, del problema
\[U = \big(U(x_0),U(x_{1/2}),\ldots,U(x_{N-1/2}),U(x_N)\big)^{\mathrm{t}},\]
y $F$ representa la restricción de $f$ a la malla  mimética:
\[F = \big(f(x_0),f(x_{1/2}),\ldots,f(x_{N-1/2}),f(x_N)\big)^{\mathrm{t}}.\]

Como el operador divergencia discretizado no actúa sobre la frontera, entonces las condiciones de contorno
Robin \eqref{ec:2} se obtienen a partir
\begin{equation}\label{ec:8}
    \big[[\alpha]+[\beta](\BB\GG)\big]U = f_b,
\end{equation}
donde el vector $f_b$ resulta de restringir el término no homogéneo de las condiciones de contorno a la
malla mimética, es decir, $f_b=(\gamma_a,0,\ldots,0,\gamma_b)^{\mathrm{t}}$. Las matrices $[\alpha]$ y $[\beta]$ son
tales que $\alpha_{1,1}=\alpha_a$,~ $\alpha_{N+2,N+2}=\alpha_b$,~ $\beta_{1,1}=\beta_a$,~ $\beta_{N+2,N+2}=\beta_b$, y el resto de entradas son cero.

A partir de \eqref{ec:7} y \eqref{ec:8}, el esquema mimético para la ecuación de convección-difusión \eqref{ec:1} sujeto a las condiciones de contorno Robin \eqref{ec:2} queda dado por
\begin{equation}\label{ec:9}
    \Big[[\alpha]+[\beta](\BB\GG)+k\DD\GG+\nu\GG\Big]U = F+f_b.
\end{equation}

\subsection{Método de Diferencias Finitas}
En el método de diferencias finitas los operadores diferenciales presentes en la ecuación diferencial son discretizados (aproximados) para los nodos de la malla, $x_i$, con $i=0,1,\ldots,N$, usando esquemas en diferencia que provienen del truncamiento del polinomio de Taylor.  Para obtener el mismo orden de convergencia que el esquema mimético \eqref{ec:7} planteado en el apartado anterior, se propone un método usando los esquemas de diferencia  de segundo orden:
\[
\dfrac{d^2}{dx^2}u(x_i) = \frac{u(x_{i+1})-2u(x_i)+u(x_{i-1})}{h^2}+\mathcal{O}(h^2)\qquad
\dfrac{d}{dx}u(x_i) = \frac{u(x_{i+1})-u(x_{i-1})}{2h}+\mathcal{O}(h^2).
\]
En los nodos interiores de la malla, la discretización realizada a partir del esquema de
segundo orden es riguroso. No obstante, al usar el esquema de diferencias centradas de segundo orden para
discretizar las condiciones de frontera de Robin \eqref{ec:2} se emplean nodos fantasmas, lo cual hace al
método poco riguroso en este respecto.

El esquema en diferencias finitas para la ecuación de convección-difusión \eqref{ec:1} queda dado entonces
\begin{eqnarray}\label{ec:10}
 \nonumber \Big[(k_0-\nu_0\frac{h}{2})\frac{2h\alpha_a}{\beta_a}-2k_0\Big]U_0+2k_0U_1 &=& h^2f_0 + (k_0-\nu_0\frac{h}{2})\frac{2h}{\beta_a}\gamma_a,\\[2mm]
  (k_i-\nu_i\frac{h}{2})U_{i-1} -2k_iU_i+(k_i+\nu_i\frac{h}{2})U_{i+1} &=&h^2f_i,\qquad\mbox{para $i=1:N-1$},\\[2mm]
  2k_NU_{N-1} - \Big[(k_N+\nu_N\frac{h}{2})\frac{2h\alpha_b}{\beta_b}+2k_N\Big]U_N &=& h^2f_N - (k_N+\nu_N\frac{h}{2})\frac{2h}{\beta_b}\gamma_b,\nonumber
\end{eqnarray}
donde $U_i$, con $i=0,1,\ldots,N$, define la solución aproximada del problema. Los subíndices en el esquema denotan la evaluación en los nodos $x_i$ de la malla de diferencias finitas (por ejemplo, $f_i=f(x_i)$).

\subsection{Método de los Elementos Finitos}
En el MEF se tiene una solución $\uu:=\sum a_i\phi_i(x)$, con $a_i$ coeficientes a determinar, y $\phi_i$ funciones base de un espacio de dimensión finita  $\VVH\subset\VV$, con $\VV$ el espacio de funciones admisibles de la forma débil o variacional del problema
\begin{equation}\label{ec:11}
   B(u,v)=l(v), \qquad\forall v\in\VV,
\end{equation}
donde
\begin{equation}\label{ec:12}
    B(u,v)=\theta_a\alpha_au(a)-\theta_b\alpha_bu(b)-\int_a^b\big[(kv)'+\nu v\big]u'dx
\end{equation}
es una forma bilineal y
\begin{equation}\label{ec:13}
   l(v)=\int_a^bfvdx-\theta_b\gamma_b+ \theta_b\gamma_b
\end{equation}
es un funcional lineal. En \eqref{ec:12} y \eqref{ec:13}, $\theta_a$ y $\theta_b$ quedan dadas por $\theta_a=k(a)v(a)/\beta_a$ y $\theta_b=k(b)v(b)/\beta_b$, y $\VV$ el espacio de Sobolev $\mathbf{H}^1(\Omega)$. Los coeficientes $a_i$ son determinados
a partir del sistema de ecuaciones que surge al sustituir $\uu$ en la forma variacional \eqref{ec:11}.

\section{Resultados Numéricos}
En este apartado, se presentan y analizan los resultados numéricos  obtenidos (comparativamente) por los  métodos propuestos en  distintas configuraciones del problemas \eqref{ec:1}-\eqref{ec:2} en el dominio $\Omega=(0,1)$. Para el análisis  se usan las normas del máximo, $\|\cdot\|_{\infty}$, y la norma $L_2$ definidas por
\[\|e\|_{\infty}:=\max\big\{e_j:=|\tilde{u}_j-u_j|\big\},\qquad\qquad \|e\|_{L_2}^2:=\int_{\Omega}|e|^2\dOmega.\]
En la norma del máximo, $\tilde{u}_j$ representa la solución obtenida por alguno de los tres métodos numéricos y $u_j$ la solución exacta del problema. Para el caso del esquema mimético, $j=i+1/2$, con $i=0,\ldots,n-1$; y $j=1,\ldots,n-1$, para los otros dos métodos numéricos. Para la norma $L_2$, se usa la forma continua en lugar de su versión discreta, $\|e\|_{L_2}^2=h\sum_{j}|\tilde{u}_j-u_j|^2$, usada por \citep{CastilloBatista2009}. Se quiere evitar el sesgo que puede presentar la versión discreta de la norma a favor del esquema mimético o diferencias finitas. Pues, en cada segmento de la malla, la solución de estos métodos es constante, mientras la de elementos finitos, no necesariamente.

\subsection*{Ejemplo 1}
\begin{figure}[!h]
 \centering
 \includegraphics[clip,angle=0,width=.4\hsize]{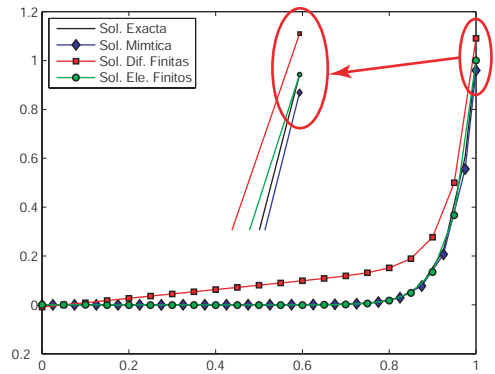}
\caption{Ejemplo 1: Comparación de las soluciones aproximadas en un malla de 20 nodos. En la frontera, el FEM
supera en precisión a los otros dos esquemas. }\label{Ejemplo1:figura1:soluciones}
\end{figure}

Se resuelve \eqref{ec:1}-\eqref{ec:2} con coeficiente de difusión $k=1$,  velocidad convectiva $\nu=0$, condiciones de contorno Robin:
\[
   \alpha u(0)-u'(0)=-20/(e^{20}-1),\qquad\qquad
    \alpha u(1)+u'(1)=0,
\]
con $\alpha:=-20e^{20}/(e^{20}-1)$ y un término fuente, $f(x)$, definido tal que la solución analítica del problema queda dada por $u(x)=(e^{20x}-1)/(e^{20}-1)$. La Figura \ref{Ejemplo1:figura1:soluciones} muestra las soluciones aproximadas junto a la solución analítica del problema en una malla uniforme de 20 puntos.
La ampliación en la frontera derecha muestra la mejor precisión del MEF al imponer condiciones de frontera,  seguido del esquema mimético que supera a DF, aunque en esta se esté usando esquemas de segundo orden para aproximar la frontera (ver \eqref{ec:10}). Esta propiedad de aproximación se mantiene al aumentar el número de nodos y se repite para cualquier problema con coeficientes $k$ y $\nu$ constantes. La igualdad en la precisión de los nodos de frontera es solo alcanzada cuando se tienen condiciones de frontera Dirichlet.

\begin{figure}[!h]
 \centering
 \includegraphics[clip,angle=0,width=.47\hsize]{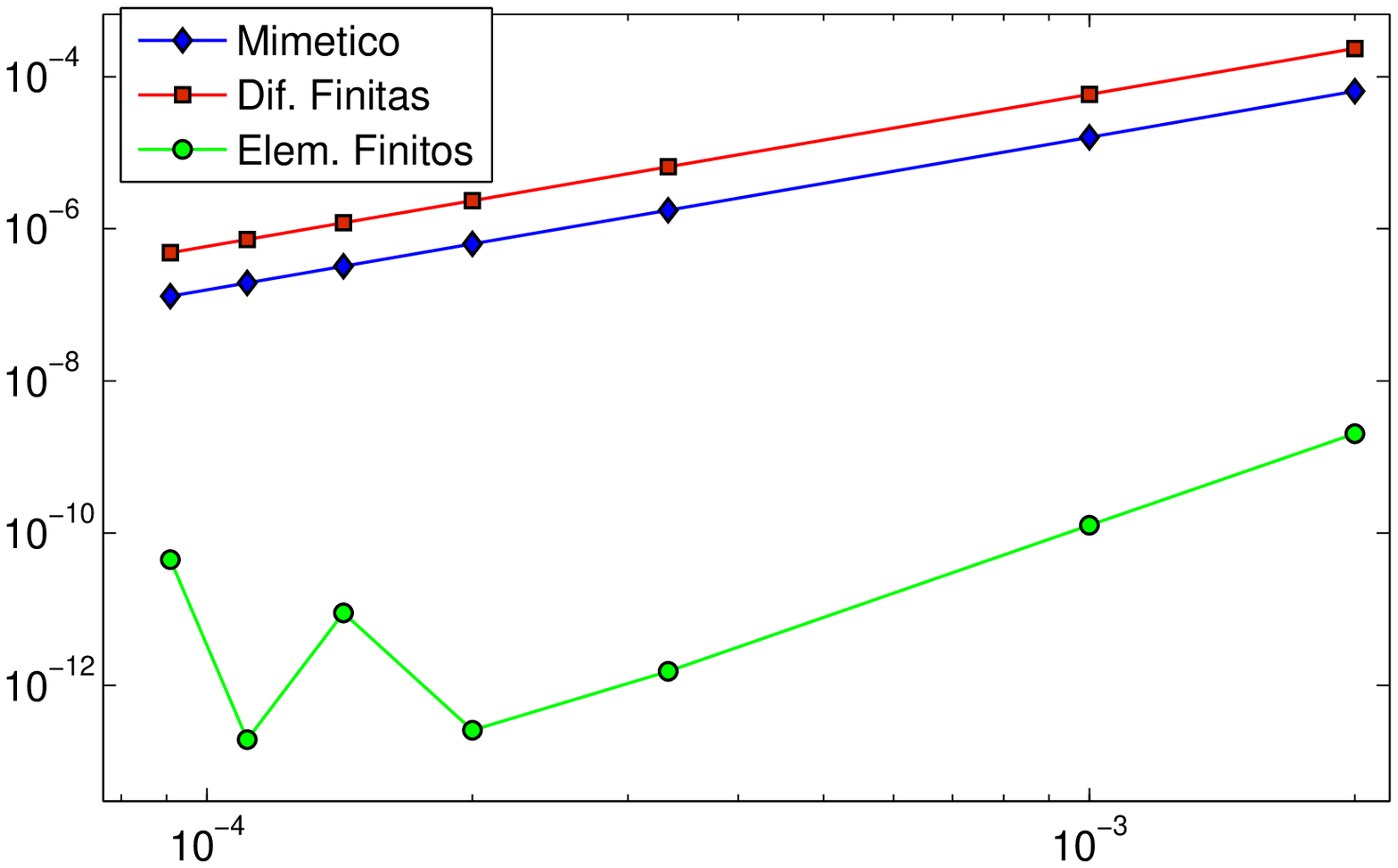}
 \includegraphics[clip,angle=0,width=.47\hsize]{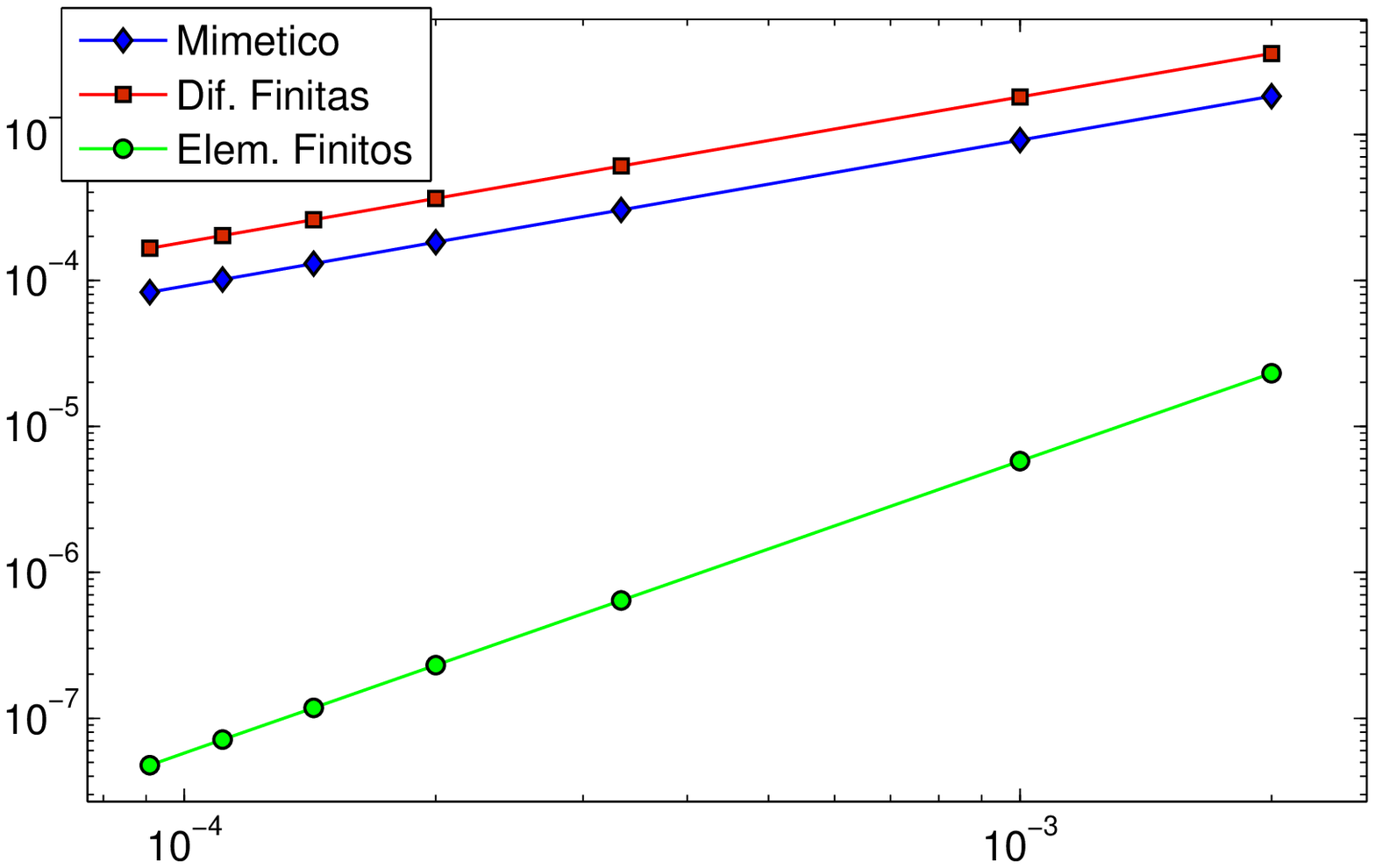}
\caption{Ejemplo 1: Errores numéricos: en noma del máximo, $\|\cdot\|_{\infty}$, (izquierda) y norma $L_2$ (derecha). }\label{Ejemplo1:figura2:Errores}
\end{figure}

La Figura \ref{Ejemplo1:figura2:Errores} muestra los errores alcanzados (a la izquierda medidos a partir de la norma del máximo y, a la derecha, usando norma $L_2$). En la gráfica, las pendientes de las rectas definen el orden de convergencia de los métodos y mientras más abajo se encuentre, mejor será el esquema que represente.
\begin{table}[!hp]
 \footnotesize
\begin{center}
 \setlength{\extrarowheight}{1.5pt}
  \begin{tabular}{|c|ccc|ccc|}\hline
  \hbox{$N^o$ de}& \multicolumn{3}{c}{Norma $\|\cdot\|_{\infty}$} &\multicolumn{3}{|c|}{Norma $\|\cdot\|_{L_2}$} \\\cline{2-7}
  \hbox{Elementos}& DF      & MIM     & MEF     & DF      & MIM     & MEF  \\\hline
  ~~1000          &1.999944 &2.016795 &4.000306 &0.986027 &1.000097 &1.999880\\
  ~~3000          &1.999990 &2.007299 &4.022705 &0.994014 &1.000020 &1.999978\\
  ~~5000          &1.999998 &2.003188 &3.485599 &0.997404 &1.000004 &1.999996\\
  ~~7000          &2.000000 &2.002136 &-10.558  &0.998307 &1.000002 &1.999967\\
  ~~9000          &1.999999 &2.001412 &15.251   &0.998740 &1.000001 &2.000043\\
  11000           &1.999998 &2.001592 &-27.111  &0.998995 &1.000001 &2.000664\\\hline
  \end{tabular}
\end{center}
\caption{Ejemplo 1. Orden de convergencia en Norma del M\'aximo y Norma $L_{2}$.}\label{Ejemplo1:Tabla1}
\end{table}
Los métodos alcanzan el orden de convergencia teórico en la norma del máximo (ver Tabla \ref{Ejemplo1:Tabla1}). Sin embargo, en norma del máximo, para este problema,  el MEF presenta una superconvergencia alcanzando el $\varepsilon$ de máquina. Por tal motivo, las oscilaciones que presenta la gráfica para el caso del MEF es debido a errores de redondeo de máquina y no a pérdida de precisión del MEF.  Esta característica no se mantiene en otras configuraciones del problema como se verá posteriormente. En norma $L_2$, el MEF presenta un orden dos en su convergencia (como era de esperar, ver \citep{Solin2006}), mientras los otros dos métodos presentan un orden uno de convergencia (Figura \ref{Ejemplo1:figura2:Errores} y Tabla \ref{Ejemplo1:Tabla1}).

\subsection*{Ejemplo 2}
\begin{figure}[!h]
 \centering
 \includegraphics[clip,angle=0,width=1.0\hsize]{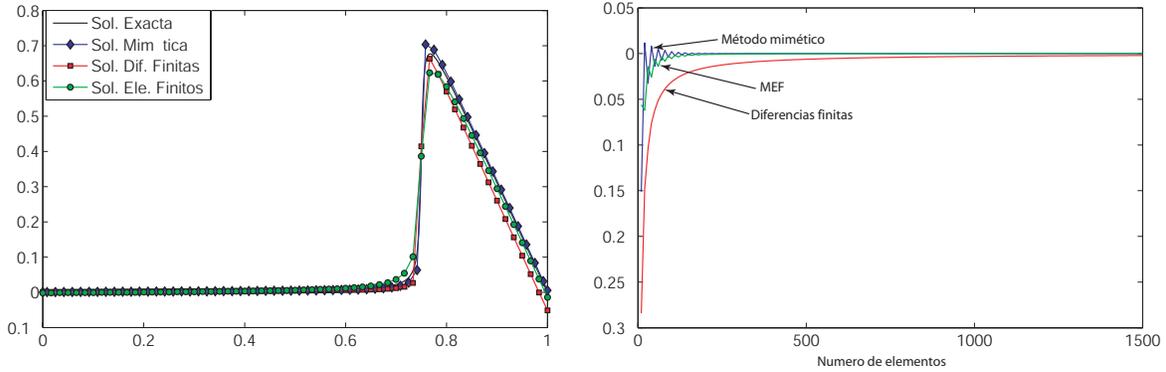}
\caption{Ejemplo 2: Soluciones aproximadas y exacta en un malla de 60 nodos (izquierda). Convergencia de los métodos en el nodo de frontera $x=1$ a medida que se aumenta el número de nodos en la malla. }\label{Ejemplo2:figura1:soluciones}
\end{figure}

En este caso, se define el problema modelo \eqref{ec:1}-\eqref{ec:2} con coeficiente de difusión $k(x)=1/\alpha+\alpha(x-x_0)^2$,  velocidad convectiva $\nu(x)=k'(x)$,  con $\alpha=250$ y $x_0=0.75$, condiciones de contorno Robin:\nopagebreak
\[
   u(0)+u'(0)=\alpha/(1+\alpha^{2}x_0^2),\qquad\qquad
   u(1)+u'(1)=-\arctan(\alpha(1-x_0))-\arctan(\alpha x_0),
\]
 y un término fuente, $f(x)$, definido tal que la solución analítica del problema viene dada por $u(x)=(1-x)\big[\arctan(\alpha(x-x_0))+\arctan(\alpha x_0)\big]$. La Figura \ref{Ejemplo2:figura1:soluciones} (izquierda) muestra las soluciones apro\-ximadas junto a la solución analítica del problema para una malla
uniforme de 60 puntos. A la derecha, se ilustra la convergencia asintótica, $h\rightarrow 0$, de las aproximaciones en el nodo de frontera $x=1$. Cuando la discretización es gruesa (pocos nodos) el método mimético supera levemente al MEF. El método DF necesita superar los 1500 elementos para alcanzar la precisión de los otros dos métodos.
Para este ejemplo, se puede ver que, $u(0)=u(1)=0$, lo cual simplifica las condiciones de frontera y, por lo tanto, la formulación variacional del problema (ver, ecuaciones \eqref{ec:12}-\eqref{ec:13}). En estos casos, la mejor precisión en los nodos de frontera es lograda por el MEF.

\begin{figure}[!h]
 \centering
 \includegraphics[clip,angle=0,width=.9\hsize]{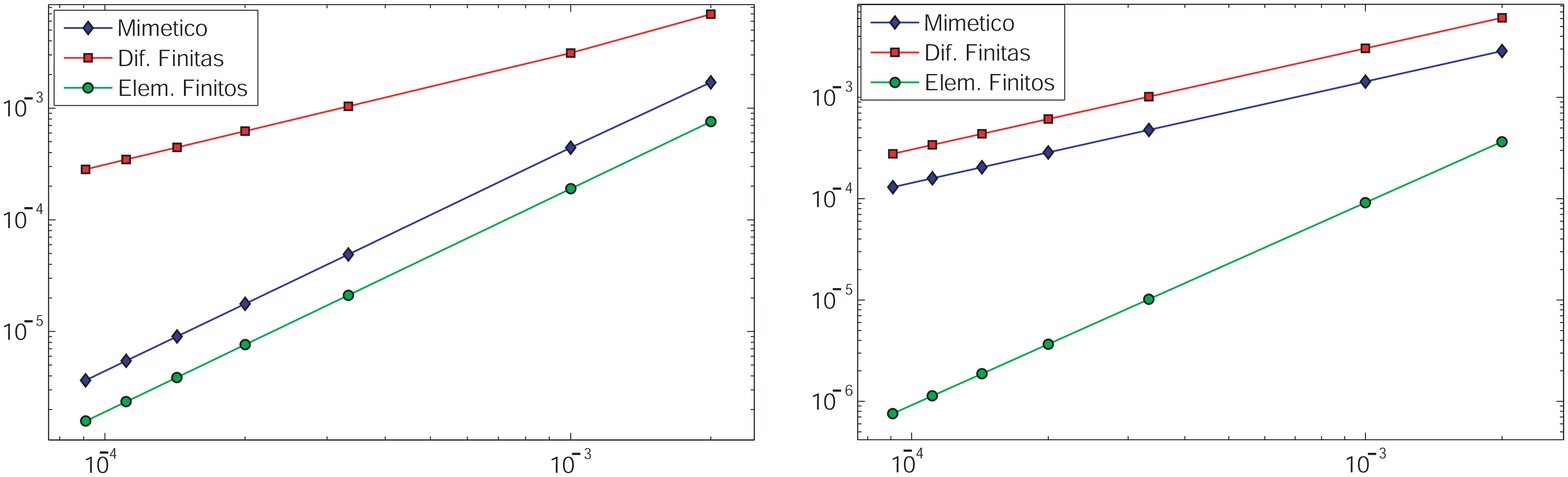}
\caption{Ejemplo 2: Errores numéricos: en noma del máximo, $\|\cdot\|_{\infty}$, (izquierda) y norma $L_2$ (derecha). }\label{Ejemplo2:figura2:Errores}
\end{figure}

La Figura \ref{Ejemplo2:figura2:Errores} y Tabla \ref{Ejemplo2:Tabla1} muestra los errores y el orden de convergencia alcanzado.
En norma del máximo, el método DF pierde el orden 2 de convergencia de los esquemas usados para su construcción. Este fenómeno es debido a que los coeficientes de difusión y convectivo dependen de la variable espacial (no constantes). El MIM y MEF mantienen el orden 2 de convergencia  con MEF, superando ligeramente en precisión a MIM. En norma $L_2$, se repite el comportamiento del ejemplo anterior: orden dos para MEF y orden uno para DF y MIM.

\begin{table}[!http]
 \footnotesize
\begin{center}
 \setlength{\extrarowheight}{1.5pt}
  \begin{tabular}{|c|ccc|ccc|}\hline
  \hbox{$N^o$ de}& \multicolumn{3}{c}{Norma $\|\cdot\|_{\infty}$} &\multicolumn{3}{|c|}{Norma $\|\cdot\|_{L_2}$} \\\cline{2-7}
  \hbox{Elementos}& DF      & MIM     & MEF     & DF      & MIM     & MEF  \\\hline
  ~~1000          &1.154080 &1.942273 &1.993329 &1.004321 &1.003528 &1.996938 \\
  ~~3000          &1.000175 &2.004354 &1.998233 &0.998419 &1.000636 &1.998873 \\
  ~~5000          &1.000073 &1.994971 &1.999533 &0.998652 &1.000109 &1.999803 \\
  ~~7000          &1.000047 &2.000917 &1.999993 &0.999021 &1.000045 &1.999919 \\
  ~~9000          &1.000035 &2.000485 &1.999783 &0.999236 &1.000025 &1.999950 \\
  11000           &1.000028 &1.998615 &2.000046 &0.999375 &1.000016 &1.999978 \\\hline
  \end{tabular}
\end{center}
\caption{Ejemplo 2. Orden de convergencia en Norma del M\'aximo y Norma $L_{2}$.}\label{Ejemplo2:Tabla1}
\end{table}

\subsection*{Ejemplo 3}
Se resuelve el problema modelo \eqref{ec:1}-\eqref{ec:2} para un coeficiente de difusión constante $k=1.052$,  velocidad convectiva $\nu=-110.5$, y un término fuente $f\equiv0$. La solución analítica queda dada por
$u(x)=(1-e^{-\lambda x})/(1-e^{-\lambda})$, con $\lambda=\nu/k$. En la experimentación numérica se toman inicialmente condiciones de contorno Dirichlet:
$u(0)=0,\quad u(1)=1$;
y, posteriormente, condiciones Robin:
\[
   u'(0)=\lambda/(1-e^{-\lambda}),\qquad\qquad
   u(1)+u'(1)=1+\lambda e^{-\lambda}/(1-e^{-\lambda})
\]

\begin{figure}[!h]
 \centering
 \includegraphics[clip,angle=0,width=1.0\hsize]{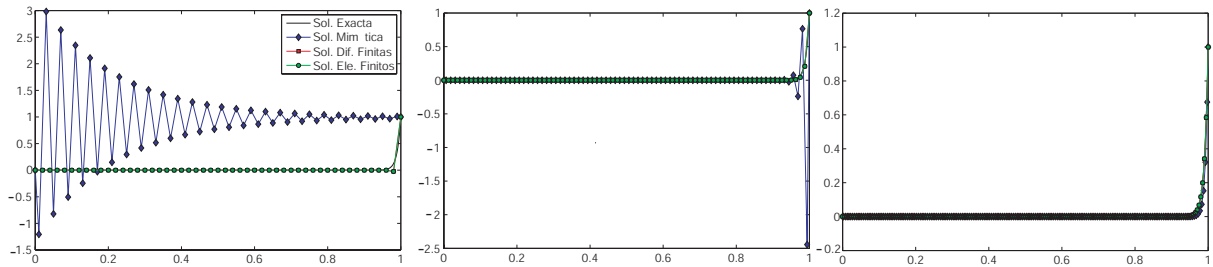}
\caption{Ejemplo 3: Soluciones aproximadas y exacta para 50 (izquierda), 80 (centro) y 200 (derecha) elementos. }\label{Ejemplo3:figura1:Soluciones_Dirichlet}
\end{figure}

Para esta configuración del problema junto a condiciones de contorno Dirichlet, es altamente conocido y fácil de probar, los métodos DF y MEF resultan equivalentes y presentan oscilaciones cuando el número local de Péclet $|\nu|h/(2k)\gg1$. La forma más simple de superar este fenómeno oscilatorio (físicamente no correcto) es  hacer $h$ suficientemente pequeño o añadiendo un término de difusividad artificial al problema. Debe quedar claro, que en los casos en que el número local de Péclet es igual o mucho menor que 1 estos métodos  no presentan oscilaciones, contrario a lo que erróneamente se afirma en \citep{CastilloBatista2009}.
La Figura \ref{Ejemplo3:figura1:Soluciones_Dirichlet} muestra las soluciones obtenidas para 50, 80 y 200 elementos en la malla. Resulta evidente que la solución del esquema mimético también presenta las oscilaciones, y para este caso necesita más nodos que DF y MEF para evitas las oscilaciones y lograr la convergencia a la solución analítica.

\begin{figure}[!h]
 \centering
 \includegraphics[clip,angle=0,width=.9\hsize]{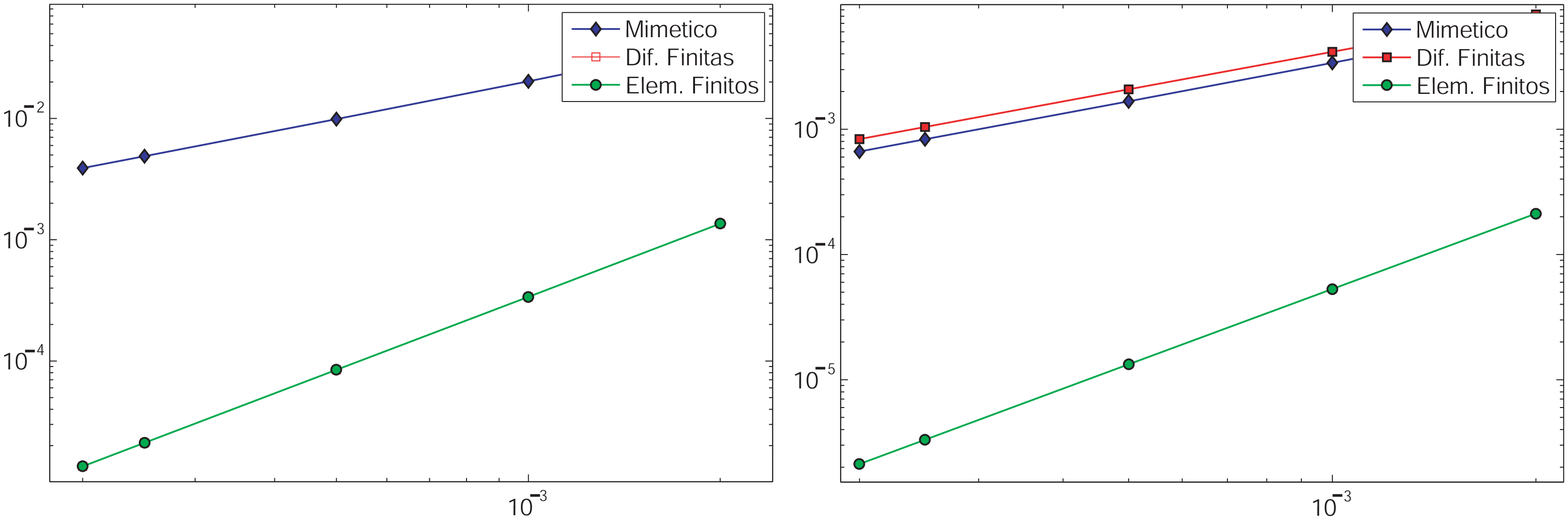}
\caption{Ejemplo 3: Errores numéricos: en noma del máximo, $\|\cdot\|_{\infty}$, (izquierda) y norma $L_2$ (derecha). }\label{Ejemplo3:figura2:ErroresDirichlet}
\end{figure}

La Figura \ref{Ejemplo3:figura2:ErroresDirichlet} y Tabla \ref{Ejemplo3:Tabla1}
muestran los errores y el orden de convergencia alcanzado para el caso de condiciones de contorno Dirichlet.
\begin{table}[!h]
 \footnotesize
\begin{center}
 \setlength{\extrarowheight}{1.5pt}
  \begin{tabular}{|c|ccc|ccc|}\hline
  \hbox{$N^o$ de}& \multicolumn{3}{c}{Norma $\|\cdot\|_{\infty}$} &\multicolumn{3}{|c|}{Norma $\|\cdot\|_{L_2}$} \\\cline{2-7}
  \hbox{Elementos}& DF      & MIM     & MEF     & DF      & MIM     & MEF  \\\hline
  ~~1000          &2.005097 &1.078854 &2.005097 &0.990024 &1.042730 &1.998561 \\
  ~~2000          &1.999493 &1.034747 &1.999493 &0.995134 &1.018963 &1.999640 \\
  ~~4000          &2.000324 &1.016534 &2.000324 &0.997599 &1.008977 &1.999910 \\
  ~~5000          &2.000267 &1.009977 &2.000267 &0.998517 &1.005461 &1.999966 \\\hline
  \end{tabular}
\end{center}
\caption{Ejemplo 3. Orden de convergencia en Norma del M\'aximo y Norma $L_{2}$.}\label{Ejemplo3:Tabla1}
\end{table}
En norma del máximo, como ya se dijo antes, los métodos DF y MEF resultan equivalentes, y mantienen el orden 2
de convergencia. Sin embargo, el esquema MIM pierde su orden 2 de convergencia y su precisión es más pobre que la lograda por los otros dos métodos. En norma $L_2$, el método DF se une a MIM logrando solo un orden uno en su convergencia, y con un error exiguo entre los dos. En esta norma, el MEF mantiene el orden dos de convergencia y nuevamente resulta superior en su precisión.
\begin{figure}[!h]
 \centering
 \includegraphics[clip,angle=0,width=1.0\hsize]{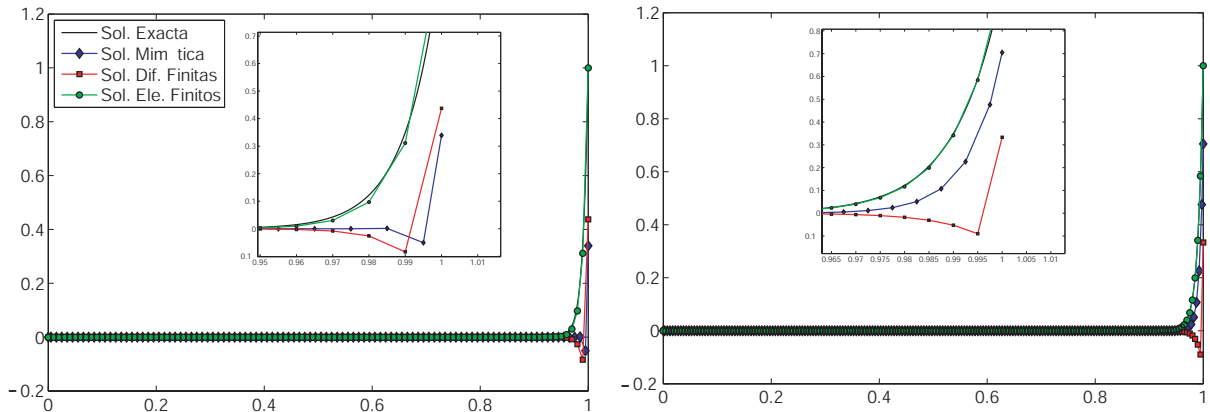}
\caption{Ejemplo 3: Soluciones aproximadas y exacta para 100 (izquierda) y 200 (derecha) elementos. La parte ampliada en cada gráfica representa el comportamiento numérico alrededor de la frontera $x=1$. }\label{Ejemplo3:figura3:Soluciones_Robin}
\end{figure}

Para las condiciones de contorno Robin, las soluciones aproximadas son mostradas en la Figura \ref{Ejemplo3:figura3:Soluciones_Robin}. En la frontera $x=1$ la precisión es pobre para los métodos DF y MIM (parte ampliada al centro de cada gráfica de la Figura \ref{Ejemplo3:figura3:Soluciones_Robin}). El esquema MIM llega a necesita hasta 4 mil elementos para alcanzar la precisión que logra MEF con solo 80 elementos. EL método DF pierde precisión en esta frontera y no converge.

Los resultados de convergencia para MEF y MIM resultan equivalentes a los alcanzados para el caso de condiciones de contorno Dirichlet (orden 2 para MEF y orden 1 para MIM). Sin embargo, diferencia finitas pierde por completo su convergencia, debido a la perdida de precisión que presenta alrededor de la frontera $x=1$.
Por último, se debe señalar que la matriz que define el sistema resultante para MEF y MIM presenta un número de condición superior a $10^{-19}$.


\section{Conclusiones y Comentarios Finales}

Resulta claro, que el MEF presenta mejores resultados en precisión y convergencia que los otros dos métodos. Al mismo tiempo, MIM resulta en todo momento superior al método DF  (como ya ha sido reportado por muchas referencias). En el caso del Ejemplo 3, donde se podría esperar superioridad del MIM, debido a su condición conservativa, tampoco logra superar a MEF. Sin querer entrar en discusión de las modificaciones apropiadas o métodos óptimos para la resolución del Ejemplo 3, queda abierta la pregunta de qué condiciones o cambios se deben imponer al MIM para problemas altamente convectivos.

En gran medida, se ha justificado el uso de MIM ante MEF debido a la complejidad teórica de MEF (para aquellos no matemáticos). Sin embargo, la definición de los operadores discretos usados por MIM (principalmente en el caso 2D o 3D) puede resultar tan complejo en cálculo que pueden competir fácilmente con la necesidad teórica de MEF. Tal dificultad se acrecienta, si se quiere subir el orden de convergencia del método, hecho que se logra en MEF con gran facilidad. Además, hay que añadir la obligatoriedad que existe de una eficiente implementación computacional, pues el  uso de matrices para definir operadores y condiciones de contorno puede resultar en un rápido desbordamiento de memoria en el método MIM.
Por último, y tal vez lo más importante, se debe mencionar el alto grado de complejidad y dificultad  que presenta MIM para trabajar en mallas que se ajusten a dominios arbitrarios o en elementos deformados \citep{Hyman1997,Hyman97}.

Como conclusión, si se quiere abordar un problema con un método novel y se tiene a mano la discretización de los operadores, el método mimético, MIM, resultará una buena elección, por arriba de cualquier esquema de diferencias finitas. Desde la vertiente matemática, sobresale el interés de afianzar las bases teóricas del método MIM y proponer mejoras en las líneas donde el mismo aún sigue sin rigor o no han sido abordadas.

\bibliography{theBiblio}
\end{document}